\documentclass{crm-article}
\usepackage{amsmath}

\usepackage{citehack}
\usepackage{algorithm}
\usepackage{algorithmic}
\usepackage{theorem}

\newtheorem{theorems}{Теорема}

\begin{document}

\year=2020 % current year by default
\journalVol{10}

%выпуска
\journalNo{1}
\setcounter{page}{1}

% раздел журнала
\journalSection{Математические основы и численные методы моделирования}
\journalSectionEn{Mathematical modeling and numerical simulation}

% дата получения
%\journalReviewed{01.06.2016.}

%принято к публикации
%\journalAccepted{01.06.2016.}

%!!!!!!!! \englishpaper раскомментировать в том случае, если текст статьи на английском языке
%\englishpaper

%!!!!!!!! \affiliationnoref раскомментировать, если автор один или все авторы из одной организации
%\affiliationnoref

%!!!!!!!! \emailnoref раскомментировать в том случае, если автор единственный
%\emailnoref

\UDC{519.85}
\title{Ускоренные адаптивные по константам сильной выпуклости и Липшица для градиента методы первого порядка}
\titleeng{Fast adaptive by constants of strong-convexity and Lipschitz for gradient first order methods}
%\thanks{Работа была поддержана грантом РФФИ 18-31-20005 мол-а-вед в первой части и грантом РНФ 17-11-01027 во второй.}%если имеется
%\thankseng{This work was supported by RFFI 18-31-20005 mol\_a\_ved in the first part of the work and by RSCF grant No. 17-11-01027 in the second part of the work.}

%автор - в формате \author{\firstname{И.\,И.}~\surname{Иванов}}
\author[1,2,3]{\firstname{Н.\,В.}~\surname{Плетнев}}
%автор - в формате \authorfull{Имя Отчество Фамилия}
\authorfull{Никита Вячеславович Плетнев}
% автор на англ. - в формате \authoreng{\firstname{I.\,I.}~\surname{Ivanov}}
\authoreng{\firstname{N.\,V.}~\surname{Pletnev}}
%автор на англ. - в формате \authorfull{Firstname M. Surname}
\authorfulleng{Nikita V. Pletnev}
%вписать свою электронную почту
\email{nikita.pletnev@phystech.edu}
%организация - в формате \affiliation{Московский государственный университет,\protect\\ Россия, 141700, г. Москва, ул. Университетская, д. 9}
\affiliation[2]{Московский физико-технический институт,\protect\\ Россия, 141701, г. Долгопрудный, Институтский пер., д. 9}
%организация - в формате \affiliationeng{Moscow State Institute University, 9 University street, Moscow, 141700, Russia}
\affiliationeng{Moscow Institute of Physics and Technology,\protect\\ 9 Institute lane, Dolgoprudny, 141701, Russia}
\affiliation[3]{Вычислительный центр им. А.А. Дородницына Российской академии наук Федерального исследовательского центра <<Информатика и управление>> Российской академии наук,\protect\\ Россия, 119333, Москва, ул. Вавилова, 40}
%организация - в формате \affiliationeng{Moscow State Institute University, 9 University street, Moscow, 141700, Russia}
\affiliationeng{Institution of Russian Academy of Sciences Dorodnicyn Computing Centre of RAS,\protect\\ Vavilov st. 40, 119333 Moscow, Russia}

% повторите блок для каждого автора;
% если авторов несколько, и автоматическая расстановка сносок от фамилий
% к организациям приводит к неправильным результатам, укажите правильный
% вариант в квадраных скобках

\begin{abstract}
Работа посвящена построению эффективных и применимых к реальным задачам методов выпуклой оптимизации первого порядка, то есть использующих только значения целевой функции и ее производных. При построении используется быстрый градиентный метод OGM-G, который является оптимальным по сложности, но при запуске требует знания констант сильной выпуклости и Липшица градиента для вычисления количества шагов и длины шага, требуемых для достижения заданной точности. Данное требование делает невозможным практическое применение указанного алгоритма. Предлагаются адаптивный по константе сильной выпуклости алгоритм ACGM, основанный на рестартах OGM-G с обновлением оценки константы сильной выпуклости, и адаптивный по константе Липшица градиента метод ALGM, в котором применение рестартов OGM-G дополнено подбором константы Липшица с проверкой условий выпуклости, используемых в методе универсального градиентного спуска. При этом устраняются недостатки исходного метода, связанные с необходимостью знания данных констант, что делает возможным практическое использование. Доказывается, что оценки сложности построенных алгоритмов являются оптимальными с точностью до постоянного множителя. Для проверки полученных результатов проводятся эксперименты на модельных функциях и реальных задачах машинного обучения.
\end{abstract}


\begin{abstracteng}
The work is devoted to the construction of efficient and applicable to real tasks first-order methods of convex optimization, that is, using only values of the target function and its derivatives. Construction uses OGM-G, fast gradient method which is optimal by complexity, but requires to know the Lipschitz constant for gradient and the strong convexity constant to determine the number of steps and step length. This requirement makes practical usage impossible. An adaptive on the constant for strong convexity algorithm ACGM is proposed, based on restarts of the OGM-G with update of the strong convexity constant estimate, and an adaptive on the Lipschitz constant for gradient ALGM, in which the use of OGM-G restarts is supplemented by the selection of the Lipschitz constant with verification of the convexity conditions used in the universal gradient descent method. This eliminates the disadvantages of the original method associated with the need to know these constants, which makes practical usage possible. Optimality of estimates for the complexity of the constructed algorithms is proved. To verify the results obtained, experiments on model functions and real tasks from machine learning are carried out.

\end{abstracteng}
\keywordeng{fast gradient method}
\keywordeng{adaptivity on the constant for strong convexity}
\keywordeng{adaptivity on the Lipschitz constant for gradient}

\maketitle

%Раздел
\paragraph{Введение}
Работа посвящена методам оптимизации первого порядка, то есть методам, использующим лишь значения функции и ее градиента.

Задачи оптимизации функций высокой размерности имеют многообразные приложения, например, в машинном обучении, управлении, экономике и энергетике. Поскольку точное решение данных задач чаще всего невозможно, необходимо применять приближённые методы. 

На сложность задачи оптимизации влияет гладкость функции, а также то, является ли она выпуклой. Как известно, для выпуклых функций локальный минимум всегда является глобальным, и необходимое условие экстремума --- равенство градиента нулю --- становится достаточным. Методы, рассматриваемые в работе, предназначены для решения задачи выпуклой оптимизации.

Методы первого порядка пользуются большой популярностью, потому что их реализация обладает относительно невысокой вычислительной сложностью: требует вычисления только значения функции, ее градиента и простейших векторных операций.

В настоящее время активно развиваются быстрые градиентные методы, основанные на следующей идее: задается число операций, строятся оптимальные для данного числа операций последовательности коэффициентов, которые используются для получения последовательности точек. Такой подход реализован в статье \cite{kim2018fessler} (метод OGM-G), а общее описание можно найти в пособии \cite{gasnikov2017universal}.

Проблема данного подхода заключается в том, что требуемое для достижения заданного результата, например, уменьшения нормы градиента вдвое, число итераций неизвестно. Поэтому для эффективного применения подобных методов необходимо оценивать это число.

В пособии \cite{gasnikov2017universal} предлагается способ оценки, но он требует знания константы сильной выпуклости $\mu$. Также там указана предложенная Ю. Е. Нестеровым в статье \cite{nesterov1989effective} идея применения быстрого градиентного метода с оцениванием данного параметра и обновлением его значения при каждом рестарте.

Предлагается применить тот же подход к оцениванию параметра $L$.

Другим недостатком данных методов является фиксированный шаг $\dfrac{1}{L}$. Его наличие приводит к замедлению сходимости, хоть и гарантирует ее наличие. И этот параметр --- константа Липшица для градиента --- неизвестен, как и константа сильной выпуклости. Идея по подбору данного параметра реализована в алгоритме универсального градиентного спуска Ю. Е. Нестерова (параграф 5 пособия \cite{gasnikov2017universal}, первоисточник --- статья \cite{nesterov2015universal}).

\paragraph{Определения и предположения}

Решается задача безусловной минимизации:
\begin{equation}
\min_{x\in \mathrm{R}^d} f(x).
\end{equation}

Предполагается, что решение
\begin{equation}
x^*=\arg\min_{x\in \mathrm{R}^d} f(x)
\end{equation}
существует, а градиент функции $f(x)$ удовлетворяет условию Липшица с константой $L>0$:
\begin{equation}
||\nabla f(x) - \nabla f(y)||\leq L||x-y||\ \forall x,y\in \mathrm{R}^d.
\end{equation}

Здесь и далее используется 2-норма: $||a||=||a||_2=\sqrt{a_1^2+\ldots+a_d^2}.$

Также считается, что функция $f(x)$ является сильно выпуклой с неизвестной нам константой $\mu>0$:
\begin{equation}
\frac{\mu}{2}||x-x^*||^2\leq f(x)-f(x^*)\leq \frac{1}{2\mu}||\nabla f(x)||^2.
\end{equation}

%здесь всё-таки \mu в знаменателе, это видно из минимума квадратичной формы в доказательстве 

Первое неравенство напрямую следует из определения, второе доказывается в соответствии с [2]:
$$f(x^*) = \min_{y} f(y)\geq \min_y \left(f(x) + \langle \nabla f(x), y-x\rangle + \dfrac{\mu}{2} ||x-y||^2\right) = f(x) - \dfrac{1}{2\mu} ||\nabla f(x)||^2.$$

В качестве невязки нахождения точки экстремума функции используется норма градиента. Критерий останова выглядит так:
\begin{equation}
||\nabla f(x)||\leq\varepsilon.
\end{equation}

Траекторией метода оптимизации называется последовательность порождаемых им точек.

Рестартом называется перезапуск метода с использованием
результата предыдущего запуска в качестве начального значения.

Константа, с которой функция удовлетворяет определению сильной выпуклости в окрестности некоторой части траектории, превосходящая $\mu$, обозначается $\mu^{loc}$. Аналогично определяется $L^{loc}$.

Алгоритм называется адаптивным по некоторому параметру, если его применение не требует никаких предположений о значении данного параметра.

\paragraph{Обзор литературы}

Использованная в работе статья \cite{kim2018fessler} посвящена построению оптимальной последовательности коэффициентов для быстрого градиентного метода. Построенный в ней метод OGM-G является оптимальным среди методов с фиксированным числом шагов, в статье доказаны оценки для его скорости сходимости и приведён пример функции, для которой улучшение этих оценок невозможно --- откуда следует оптимальность. Изложению результатов \cite{kim2018fessler} в части оценок, касающихся целей работы, посвящен раздел 5. 

В пособии \cite{gasnikov2017universal} излагается современное состояние быстрых градиентных методов. Среди прочего, в нем содержатся идеи адаптивного подбора неизвестных констант; на этих идеях построено основное содержание работы --- разделы 6 и 7. В частности, в параграфе 5 пособия изложен придуманный Ю. Е. Нестеровым (первоисточник \cite{nesterov2015universal}) метод адаптивного подбора константы Липшица, известный как универсальный градиентный спуск. На основе данного метода в работе построен быстрый алгоритм, адаптивный как по константе сильной выпуклости, так и по константе Липшица для градиента.

В статье \cite{gasnikov2018entropy} разрабатываются численные методы оптимизации энтропии. Там применен подход, похожий на использованный в работе для подбора константы сильной выпуклости, в котором подбираемый параметр изменяется в одно и то же число раз $\beta$, и суммарное количество шагов оценивается с помощью суммы геометрической прогрессии. После этого из соображений минимизации данной оценки выбирается $\beta$.

Пособие \cite{nesterov2010introductory} посвящено изложению основ теории оптимизации и методов, которые служат основой для современных алгоритмов. Статья \cite{nesterov1989effective} содержит используемую в работе идею подбора оценок параметров.

Статья \cite{lei2019jordan} посвящена решению близкой задачи --- построению адаптивного по константе сильной выпуклости метода стохастического градиентного спуска. Однако в ней цель полностью не достигнута, так как полученный алгоритм является эффективным лишь с точностью до логарифмического множителя.

В статье \cite{fercoq2016qu} также строится метод, адаптивный по константе сильной выпуклости, но оценка сложности полученного алгоритма тоже содержит логарифмический множитель.

Работа \cite{barre2020taylor} предлагает метод, основанный на рестартах с использованием невязки по функции и адаптивный по константе сильной выпуклости. Соответственно, для его применения требуется знание целевого значения функции.

\paragraph{Исходный алгоритм}
\subparagraph{OGM-G}

В качестве базового алгоритма взят ускоренный градиентный метод с фиксированным шагом OGM-G. В \cite{kim2018fessler} показана его оптимальность в классе методов с заданным числом шагов фиксированной длины.

\begin{algorithm}
\caption{\bf{Optimal Gradient Method OGM-G}}
\label{OGM-G}
\hspace*{\algorithmicindent} \textbf{Input: }$f\in \mathcal{F}_L(\mathbb{R}^d)$, ${\bf x}_0\in\mathbb{R}^d$~--- начальная точка, $N\geq 1$.
\begin{algorithmic}[1]
\STATE ${\bf y}_0 := {\bf x}_0$
\FOR{$k=0\ldots N-1$}
\STATE ${\bf y}_{i+1} := {\bf x}_i - \frac{1}{L}\nabla f({\bf x}_i)$;
\STATE ${\bf x}_{i+1} := {\bf y}_{i+1} + \beta_i({\bf y}_{i+1} - {\bf y}_i) + \gamma_i({\bf y}_{i+1} - {\bf x}_i).$
\ENDFOR
\end{algorithmic}
\hspace*{\algorithmicindent} \textbf{Output: } $x_N$.
\end{algorithm}

Коэффициенты $\beta_i, \gamma_i$ вычисляются по формулам:
$$\beta_i = \dfrac{(\theta_i-1)(2\theta_{i+1}-1)}{\theta_i(2\theta_i-1)}; \gamma_i = \dfrac{2\theta_{i+1}-1}{2\theta_i-1},$$
где последовательность $\{\theta_i\}_{i=0}^N$ строится следующим образом:
\[ \theta_i = \begin{cases}
\frac{1+\sqrt{1+8\theta_1^2}}{2}, & i=0;\\
\frac{1+\sqrt{1+4\theta_{i+1}^2}}{2}, & 1\leq i<N;\\
1, & i=N.
\end{cases} \]

\subparagraph{Оценки}

По теореме 2 из [1], при применении OGM-G
\begin{equation}
||\nabla f(x^N)||^2\leq \frac{4L(f(x^0)-f(x^*))}{N^2}.
\end{equation}

Оттуда же,
\begin{equation}
f(x^N)-f(x^*)\leq \frac{L||x^0-x^*||^2}{N^2}.
\end{equation}

Из (4) и (6):
\begin{equation}
||\nabla f(x^N)||^2\leq \frac{4L}{N^2}\frac{1}{2\mu}||\nabla f(x^0)||^2,
\end{equation}
или
\begin{equation}
||\nabla f(x^N)||\leq \sqrt{\frac{2L}{\mu N^2}}||\nabla f(x^0)||.
\end{equation}

Таким образом, выполнение $N$ итераций гарантирует уменьшение нормы градиента $f(x)$ как минимум вдвое, где
\begin{equation}
    N=2\sqrt{2\frac{L}{\mu}}
\end{equation}

Согласно лемме 4 статьи \cite{kim2018fessler}, полученная оценка является неулучшаемой в худшем случае.

\subparagraph{Недостатки метода}

Полученная оценка показывает, что использование OGM-G неявно предполагает, помимо наличия известной константы Липшица, знание константы сильной выпуклости.

Практически во всех реальных случаях применения методов оптимизации ни одно из этих предположений не выполняется: свойства функции заранее неизвестны, а вычисление данных параметров требует нахождения минимума и максимума собственных значений матрицы Гессе, что может быть даже сложнее, чем исходная задача оптимизации.

Указанные соображения делают оптимальный теоретически метод неприменимым на практике. Решению данной проблемы посвящен следующий раздел.

\paragraph{Адаптивность по константе сильной выпуклости}

\subparagraph{ACGM}

Ю. Е. Нестеровым в работе \cite{nesterov1989effective} предложен способ построения адаптивного по $\mu$ алгоритма, основанного на рестартах.

В этом разделе OGM-G используется в качестве <<черного ящика>>, получающего на вход функцию $f$, начальную точку ${\bf x}_0$, константу Липшица $L$ и затравочную (начальную) константу сильной выпуклости $\mu_0$. Число итераций $N$, используемое методом, вычисляется по формуле (10). В дальнейшем применение OGM-G как шага в алгоритмах будет обозначаться как $OGMG(f, {\bf x}_0, L, \mu_0)$.

ACGM --- Adaptive by constant of strong Convexity Gradient Method --- решает проблему неизвестности $\mu$, инициализируя ее произвольным значением с последующим изменением.

На каждом шаге предполагаемое значение $\mu$ умножается на одно и то же $\beta>1$.

\begin{algorithm}
\caption{\bf{Adaptive by strong Convexity Gradient Method ACGM}}
\label{ACGM}
\hspace*{\algorithmicindent} \textbf{Input: } $f\in \mathcal{F}_L(\mathbb{R}^d)$, ${\bf x}_0\in\mathbb{R}^d$ --- начальная точка, $L$, $\mu_0, \beta$, $\varepsilon$.
\begin{algorithmic}[1]
\FOR{$k\geq 0$}
\IF {$||\nabla f(x_k)||\leq \varepsilon$}
\STATE $\textbf{break}$;
\ENDIF
\STATE $\mu_k:=\beta\mu_{k-1}$;
\STATE ${\bf x}_k = OGMG(f, {\bf x}_{k-1}, L, \mu_k)$;
\IF {$||\nabla f(x_k)||\leq \frac{1}{2}||\nabla f(x_{k-1})||$}
\STATE $\textbf{continue}$;
\ENDIF
\STATE $\mu_k:=\frac{\mu_k}{\beta}$
\IF {$||\nabla f(x_k)|| < ||\nabla f(x_{k-1})||$}
\STATE ${\bf x}_{k-1}\leftarrow{\bf x}_k$;
\ENDIF
\STATE $\textbf{goto 6}$:
\ENDFOR
\end{algorithmic}
\hspace*{\algorithmicindent} \textbf{Output: } $[x_0\ldots x_N]$.
\end{algorithm}

\subparagraph{Оценки}

В результате применения ACGM очередное уменьшение вдвое нормы градиента будет выполнено за
$$2\sqrt{2\frac{L}{\mu_k^{init}}}+2\sqrt{2\frac{L}{\mu_k^{init}/\beta}}+\ldots+2\sqrt{2\frac{L}{\mu_k^{init}/\beta^m}} = \sqrt{8\frac{L}{\mu_k}}\sum\limits_{i=0}^m\frac{1}{\sqrt{\beta^i}} \lesssim \frac{\sqrt{8\beta}}{\sqrt{\beta}-1}\sqrt{\frac{L}{\mu_k}}$$
итераций метода OGM-G, где $m$ --- количество повторений цикла на шаге $k$, а индекс $init$ указывает на то, что в формуле используется не конечное значение переменной, а то, которым она была инициализирована.

При этом $\mu_k$ отличается не более чем в $\beta$ раз от $\mu^{loc}$. Использование значения $\mu$, подходящего для всего пространства, могло бы повысить количество операций.

Действительно, если последовательные $s$ точек траектории ACGM лежат в области, в которой $f({\bf x})$ сильно выпукла с константой $\mu^{loc}\geq \beta^s\mu$, то в данных точках ACGM применяется с $\mu_0\geq \dfrac{\mu^{loc}}{\beta}\geq \beta^{s-1}\mu$. Тогда количество обращений к вычислению градиента для каждого уменьшения его нормы вдвое оказывается не более $2\sqrt{2\frac{L}{\beta^{s-1}\mu}}$ --- то есть, в $\beta^{\frac{s-1}{2}}$ раз меньше, чем при $\mu_k\equiv\mu$.

Суммарное количество итераций при работе ACGM с использованием критерия останова (5) оценивается следующим образом. Требуется выполнить $\\ K=\log_2\frac{||\nabla f(x^0)||}{\varepsilon}$ шагов. Каждый шаг содержит $O\left(\sqrt{\frac{L}{\mu}}\right)$ итераций, поэтому алгоритм завершит работу, выполнив $O\left(\sqrt{\frac{L}{\mu}}\log_2\frac{||\nabla f(x^0)||}{\varepsilon}\right)$ итераций --- то есть, вычислений $f(x)$ и $\nabla f(x)$.

Как показано выше, полученная оценка по порядку величины может быть уточнена. Если $\mu_k \geq \frac{\mu_k^{loc}}{\beta}$, то каждый шаг ACGM содержит не более $\frac{\sqrt{8\beta}}{\sqrt{\beta}-1}\sqrt{\frac{L}{\mu_k}}\leq \frac{\sqrt{8}\beta}{\sqrt{\beta}-1}\sqrt{\frac{L}{\mu_k^{loc}}}$ итераций, соответственно общее количество итераций не превосходит
$$\frac{\sqrt{8}\beta}{\sqrt{\beta}-1}\sum\limits_{k=0}^{K-1} \sqrt{\frac{L}{\mu_k^{loc}}}.$$

Данное вычисление основано на идее оценки из \cite{gasnikov2018entropy}.

Минимизация зависящего от $\beta$ коэффициента дает $\beta=4$.

Это значение является оптимальным лишь с точки зрения худшего случая, когда $\mu_k=\frac{\mu_k^{loc}}{\beta}$. В реальных случаях, поскольку данное равенство является лишь теоретически возможным предельным случаем, значение коэффициента может оказаться меньше данного, но в любом случае оно превосходит $\inf_{\beta>1}\frac{\sqrt{\beta}}{\sqrt{\beta}-1}=1$

Таким образом, доказаны теоремы о сходимости построенного метода.

\begin{theorems}
Алгоритм ACGM с оптимальным $\beta=4$ и $\mu_0>\max_k \mu_k^{loc}$ достигает точки ${\bf x}$, удовлетворяющей критерию останова (5), за не более чем $C\sum\limits_{k=0}^{K-1} \sqrt{\dfrac{L}{\mu_k^{loc}}}$ вычислений градиента, где $K=\log_2\frac{||\nabla f({\bf x}^0)||}{\varepsilon}$, $C=\frac{\sqrt{8}\beta}{\sqrt{\beta}-1}=8\sqrt{2}$.
\end{theorems}

Данная теорема имеет лишь теоретический смысл, поскольку оценка $\mu_k^{loc}$ крайне затруднительна. Следующая теорема содержит менее точную, но более удобно применимую оценку.

Так как $\mu_k^{loc}\geq \mu$ при всех $k$, каждое слагаемое в сумме из теоремы 1 не превосходит $\sqrt{\frac{L}{\mu}}$, откуда сразу следует

\begin{theorems}
Алгоритм ACGM с оптимальным $\beta=4$ и $\mu_0>\max_k \mu^{loc}$ достигает точки ${\bf x}$, удовлетворяющей критерию останова (5), за не более чем $CK \sqrt{\dfrac{L}{\mu}}$ вычислений градиента, где $K=\log_2\frac{||\nabla f({\bf x}^0)||}{\varepsilon}$, $C=\frac{\sqrt{8}\beta}{\sqrt{\beta}-1}=8\sqrt{2}$.
\end{theorems}

\subparagraph{К выбору оптимального $\mu_0$. Случай $\mu_0<\mu^{loc}$}

Для упрощения вычислений, пусть $\mu_0 = \frac{\mu^{loc}}{4^k}$. Тогда по формуле (10) на первом шаге ACGM будет выполнено $M=2\sqrt{2\frac{L\cdot 4^k}{\mu^{loc}}} = 2^kN$ итераций. При этом, согласно (9), норма градиента умножится не более, чем на $\sqrt{\frac{2L}{\mu^{loc} M^2}} = \frac{1}{2^{k+1}}$.

Для достижения такого результата требуется $k+1$ рестартов, то есть $(k+1)N$ итераций при использовании OGM-G.

При использовании ACGM $i$-ый рестарт выполняется с $\mu = \mu_0\beta^i$, т. е. потребует $\frac{N}{2^i}$ итераций. Суммарное количество не превосходит $2N$.

Поскольку $2^k>2$, применение заниженного значения константы сильной выпуклости приводит к значительному увеличению количества итераций.

Объединяя все сказанное о величине начального предполагаемого значения константы сильной выпуклости, оптимальным будет такой выбор $\mu_0$, что $\mu_0>\mu_k^{loc}$ при всех $k$. Таким является, например $\mu_0=L$.

\subparagraph{Иллюстрации}

Работа алгоритма иллюстрируется на двух функциях: штрафная функция логистической регрессии с регуляризацией, которая возникает из приложений и свойства которой будут подробнее рассмотрены в разделе <<Эксперименты>>, и квадратичная форма с плохой обусловленностью и потому подходящая для проверки градиентных методов. Для каждой функции приведены график зависимости логарифма нормы градиента от номера итерации и траектория метода.

\begin{figure}[h]
	\begin{minipage}[h]{0.49\linewidth}
		\center{\includegraphics[width=0.95\linewidth]{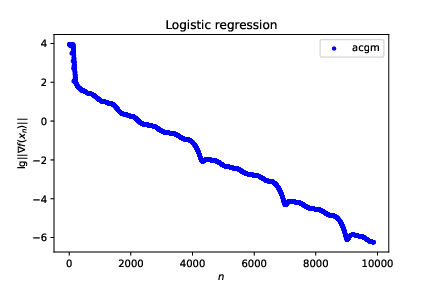} \\ Убывание нормы градиента}
	\end{minipage}
	\hfill
	\begin{minipage}[h]{0.49\linewidth}
		\center{\includegraphics[width=0.95\linewidth]{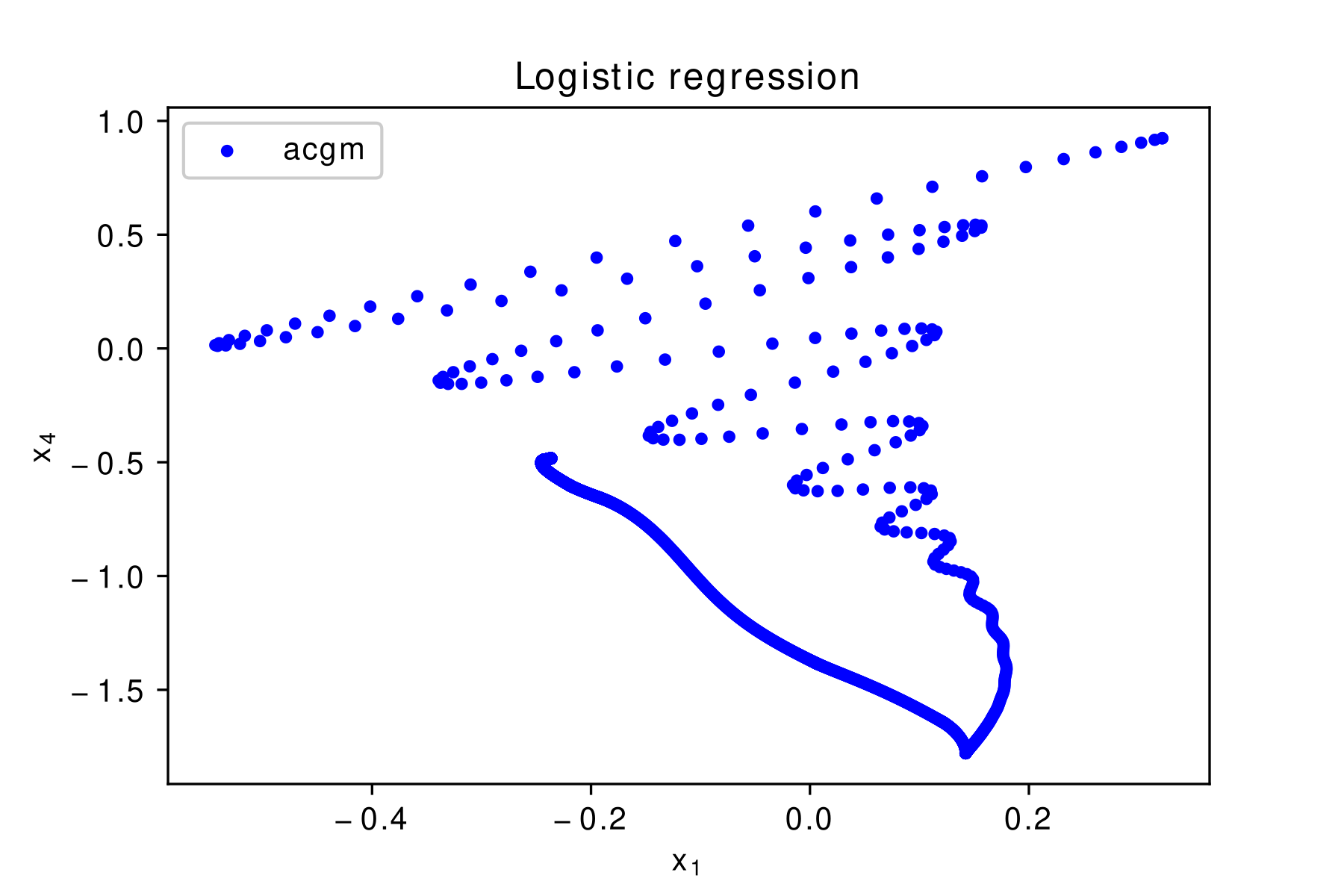} \\ Траектория метода}
	\end{minipage}
	\caption{Работа ACGM на функции логистической регрессии.}
	\label{ris:image1}
\end{figure}

\begin{figure}[h]
	\begin{minipage}[h]{0.49\linewidth}
		\center{\includegraphics[width=0.95\linewidth]{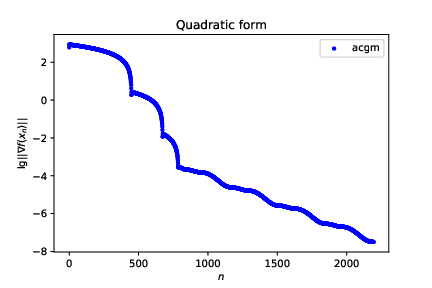} \\ Убывание нормы градиента}
	\end{minipage}
	\hfill
	\begin{minipage}[h]{0.49\linewidth}
		\center{\includegraphics[width=0.95\linewidth]{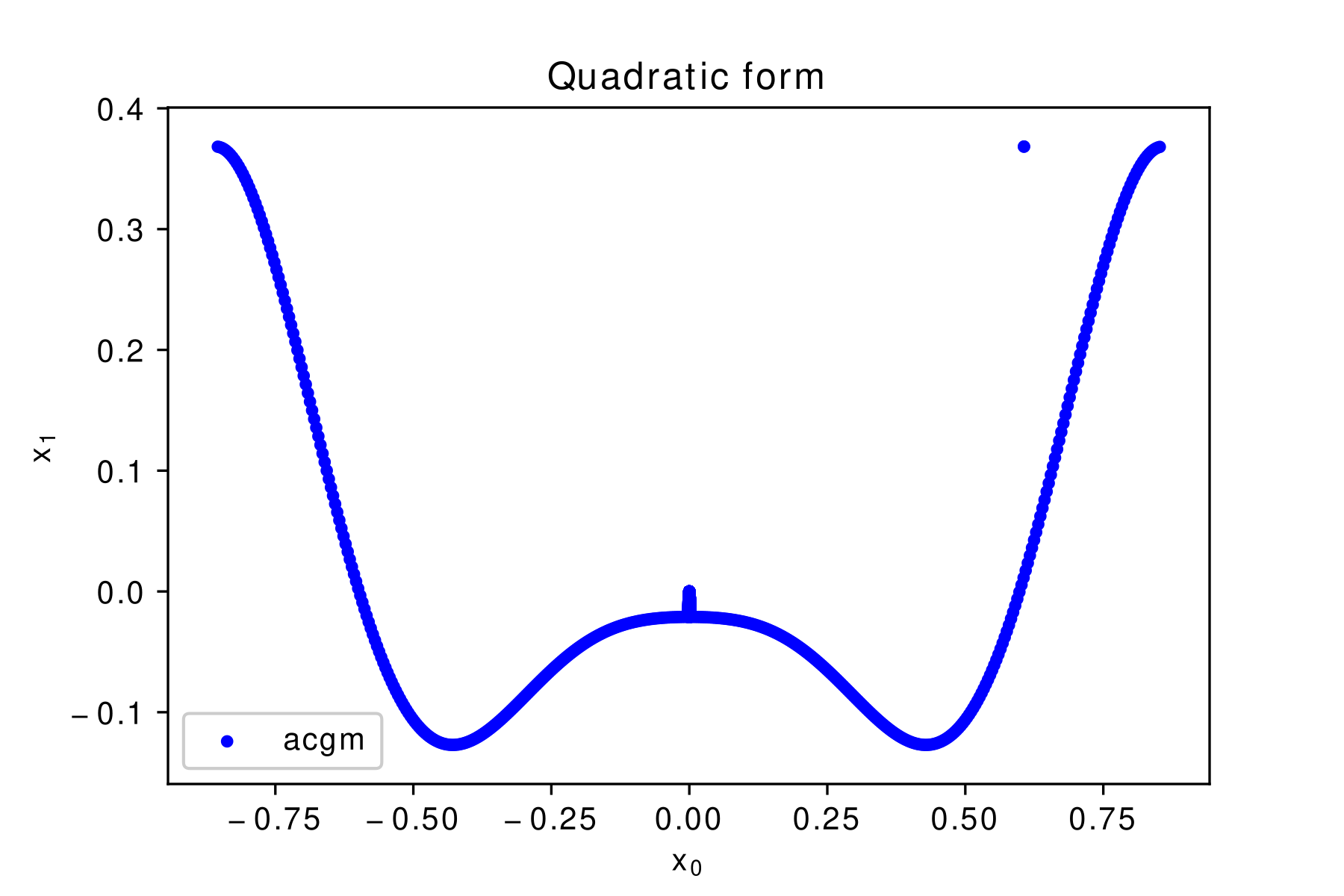} \\ Траектория метода}
	\end{minipage}
	\caption{Работа ACGM на квадратичной функции.}
	\label{ris:image2}
\end{figure}

Графики показывают, что метод сходится --- логарифм невязки убывает приблизительно линейно по количеству итераций.

\paragraph{Адаптивность по константе Липшица}
\subparagraph{Универсальный градиентный спуск}

\begin{algorithm}
\caption{\bf{Universal Gradient Method UGM}}
\label{UGM}
\hspace*{\algorithmicindent} \textbf{Input: } $f\in \mathcal{F}(\mathcal{Q})$, ${\bf x}_0\in\mathcal{Q}$, $L_0$, $\varepsilon$.
\begin{algorithmic}[1]
\FOR{$k\geq 0$}
\IF {$||\nabla f(x_k)||\leq \varepsilon$}
\STATE $\textbf{break}$;
\ENDIF
\STATE $L_{k+1}:=\frac{L_k}{2}$;
\STATE ${\bf x}_{k+1} = \arg\min_{{\bf x}\in \mathcal{Q}} \left\{f({\bf x}_k) + \langle\nabla f({\bf x}_k), {\bf x}-{\bf x}_k\rangle + L_{k+1} V({\bf x}, {\bf x}_k) \right\}$;
\IF {$f({\bf x}_{k+1}) \leq f({\bf x}_k) + \langle\nabla f({\bf x}_k), {\bf x}_{k+1}-{\bf x}_k\rangle + L_{k+1} V({\bf x}_{k+1}, {\bf x}_k)$}
\STATE $\textbf{continue}$;
\ELSE
\STATE $L_{k+1} := 2L_{k+1}$;
\STATE $\textbf{goto 6}$;
\ENDIF
\ENDFOR
\end{algorithmic}
\hspace*{\algorithmicindent} \textbf{Output: } $[x_0\ldots x_N]$.
\end{algorithm}

В замечании 2.1 пособия \cite{gasnikov2017universal} показано, что в качестве $V(x,y)$ подходит функция $V(x,y) = \frac{1}{2}||x-y||^2$.

Поскольку рассматривается задача безусловной оптимизации, $\mathcal{Q} = \mathbb{R}^d$. Градиент минимизируемого выражения равен $\nabla f({\bf x}_k) + L_{k+1}({\bf x} - {\bf x}_k)$, поэтому формула шага 6 преобразуется к виду ${\bf x}_{k+1} = {\bf x}_k - \frac{1}{L_{k+1}}\nabla f({\bf x}_k)$, а условие перехода к следующему шагу --- к виду $f({\bf x}_{k+1}) \leq f({\bf x}_k) - \frac{1}{2L_{k+1}}||\nabla f({\bf x}_k)||^2$.

\subparagraph{OGM-GL}

Универсальный градиентный спуск решает проблему неизвестности константы Липшица, но обладает недостатками простейшего градиентного метода, так как для функций с большим числом обусловленности матрицы Гессе направление градиента значительно отличается от направления на экстремум. Поэтому применение универсального градиентного метода на практике неэффективно.

Предлагается следующий вариант OGM-G, адаптивный по константе Липшица, основанный на подборе $L$ аналогично тому, как это делается в универсальном градиентном методе, с проверкой условия строки 7 при каждом вычислении ${\bf y}_{i+1}$.

Если условие нарушено, то вычисление последовательностей ${\bf x}_i$ и ${\bf y}_i$ начинается сначала с тем же количеством шагов и увеличенным значением $L$.

\begin{algorithm}
\caption{\bf{OGM-G Lipschitz}}
\label{OGM-GL}
\hspace*{\algorithmicindent} \textbf{Input: }$f\in \mathcal{F}_L(\mathbb{R}^d)$, ${\bf x}_0\in\mathbb{R}^d$~--- начальная точка, $N\geq 1$.
\begin{algorithmic}[1]
\STATE $L := \frac{L}{2}$;
\STATE ${\bf y}_0 := {\bf x}_0$;
\FOR{$i=0\ldots N-1$}
\STATE ${\bf y}_{i+1} := {\bf x}_i - \frac{1}{L}\nabla f({\bf x}_i)$;
\IF {$f({\bf y}_{i+1}) > f({\bf x}_i) - \frac{1}{2L}||\nabla f({\bf x}_i)||^2$}
\STATE $L := 2L$
\STATE $\textbf{goto 3}$
\ENDIF
\STATE ${\bf x}_{i+1} := {\bf y}_{i+1} + \beta_i({\bf y}_{i+1} - {\bf y}_i) + \gamma_i({\bf y}_{i+1} - {\bf x}_i).$
\ENDFOR
\end{algorithmic}
\hspace*{\algorithmicindent} \textbf{Output: } $x_N, L_{end}$.
\end{algorithm}

%да, goto 3, поскольку цикл начинается сначала

Коэффициенты $\beta_i, \gamma_i$ вычисляются по тем же формулам, что и в алгоритме 1. 

\subparagraph{ALGM}

\begin{algorithm}
\caption{\bf{Adaptive by Lipschitz constant Gradient Method ALGM}}
\label{ALGM}
\hspace*{\algorithmicindent} \textbf{Input: } $f\in \mathcal{F}_L(\mathbb{R}^d)$, ${\bf x}_0\in\mathbb{R}^d$ --- начальная точка, $L_0$, $\mu_0, \beta$, $\varepsilon$.
\begin{algorithmic}[1]
\FOR{$k\geq 0$}
\IF {$||\nabla f(x_k)||\leq \varepsilon$}
\STATE $\textbf{break}$;
\ENDIF
\STATE $\mu_k:=\beta\mu_{k-1}$;
\STATE ${\bf x}_k, L_k = OGMGL\left(f, {\bf x}_{k-1}, L_{k-1}, \left\lceil \sqrt{\dfrac{8L_{k-1}}{\mu_k}}\right\rceil, \varepsilon, \right)$;
\STATE $\mu_k := \mu_k\cdot\frac{L_{k_{}}}{L_{k-1}}$;
\IF {$||\nabla f(x_k)||\leq \frac{1}{2}||\nabla f(x_{k-1})||$}
\STATE $\textbf{continue}$;
\ENDIF
\STATE $\mu_k:=\frac{\mu_k}{\beta}$
\IF {$||\nabla f(x_k)|| < ||\nabla f(x_{k-1})||$}
\STATE ${\bf x}_{k-1}\leftarrow{\bf x}_k$;
\ENDIF
\STATE $\textbf{goto 6}$:
\ENDFOR
\end{algorithmic}
\hspace*{\algorithmicindent} \textbf{Output: } $[x_0\ldots x_N]$.
\end{algorithm}

Алгоритм построен на том же принципе, что и ACGM, только вместо OGM-G используется адаптивный по $L$ OGM-GL. $\mu$ изменяется для сохранения отношения $\frac{L}{\mu}$ и $N$.

\subparagraph{Оценки}

Алгоритм подбора $L$ гарантирует, согласно комментарию к алгоритму универсального градиентного спуска в пособии \cite{gasnikov2017universal}, выполнение условия $L_{k+1}\geq \dfrac{L^{loc}}{2}$, то есть $L^{loc}\leq 2L_{k+1}$. Поскольку при $L' < L''$ $\mathcal{F}_{L'}(\mathbb{R}^d) \subset \mathcal{F}_{L''}(\mathbb{R}^d)$, для OGM-GL выполнены оценки сходимости, доказанные для OGM-G, откуда следует возможность применения OGM-GL как составной части для алгоритма, подобного ACGM.

Один запуск OGM-GL требует не более $2N$ вычислений функции и $N$ вычислений градиента на каждое увеличение $L_k$, а суммарное количество вычислений функции и градиента за один запуск составляет $O\left(N\left(\log_2\dfrac{L_{end}}{L_{init}/2}+1\right)\right)$ , то есть $O\left(\sqrt{\frac{L}{\mu}}\log_{2}\dfrac{4L_k^{end}}{L_k^{init}}\right)$, т. к. $L_{loc}\leq L$, $\mu_{loc}\geq L$.

Числовой множитель составляет $\sqrt{8}$ для количества вычислений градиента и $2\sqrt{8}$ для количества вычислений значения функции.

Верхняя оценка количества итераций для каждого уменьшения нормы градиента вдвое определяется аналогично вычислению из доказательства теоремы о сходимости ACGM (расчет для количества обращений к градиенту функции; $j$ --- номер перезапуска OGM-GL в пределах одного шага ALGM):
$$\sum\limits_{j=0}^J\sqrt{8\frac{L}{\mu_k^{init}/\beta^j}}\left(2+\log_2\dfrac{L_{kj}}{L_{k,j-1}}\right)\leq \sqrt{8\frac{L}{\mu_k}} \sum\limits_{j=0}^J\frac{1}{\sqrt{\beta^j}}\left(2+\left(\log_2\dfrac{L_k}{L_{k-1}}\right)_+\right)\lesssim$$
$$\lesssim \dfrac{\sqrt{8\beta}}{\sqrt{\beta} - 1}\sqrt{\dfrac{L}{\mu_k}}\left(2+\left(\log_2\dfrac{L_k}{L_{k-1}}\right)_+\right)\leq \dfrac{\sqrt{8} \beta}{\sqrt{\beta} - 1} \sqrt{\dfrac{L}{\mu}}\left(2+\left(\log_2\dfrac{L_k}{L_{k-1}}\right)_+\right)$$

Положительная срезка логарифма появляется для того, чтобы оценка выполнялась даже в том случае, если $L_k=\frac{L_{k-1}}{2}$. Алгоритм построен так, что за запуск OGM-GL константа Липшица может уменьшиться не более, чем в два раза.

Количество шагов ALGM не превышает $K=\log_2\frac{||\nabla f({\bf x}^0)||}{\varepsilon}$. Поскольку для $t\geq -1$ выполнено свойство $(t)_+\leq t+1$, суммарное количество вычислений градиента не превосходит $\dfrac{\sqrt{8} \beta}{\sqrt{\beta} - 1} \sqrt{\dfrac{L}{\mu}}\left(3K+\log_2\dfrac{L}{L_0}\right)$; оценка количества вычислений функции отличается только числовым множителем и превышает полученное значение вдвое.

Как и для ACGM, числовой множитель минимален при $\beta=4$.
Таким образом, получена
\begin{theorems}
Траектория ALGM (алгоритма 5) с оптимальным $\beta=4$ содержит точку ${\bf x}$, удовлетворяющую критерию останова (5), после выполнения не более чем $C \sqrt{\dfrac{L}{\mu}} \left(3K+\log_2\dfrac{L}{L_0}\right)$ вычислений градиента и $2C \sqrt{\dfrac{L}{\mu}} \left(3K+\log_2\dfrac{L}{L_0}\right)$ вычислений функции, где $K=\log_2\frac{||\nabla f({\bf x}^0)||}{\varepsilon}$, $C=\frac{\sqrt{8}\beta}{\sqrt{\beta}-1}=8\sqrt{2}$.
\end{theorems}

\subparagraph{Иллюстрации}

\begin{figure}[h]
	\begin{minipage}[h]{0.49\linewidth}
		\center{\includegraphics[width=0.95\linewidth]{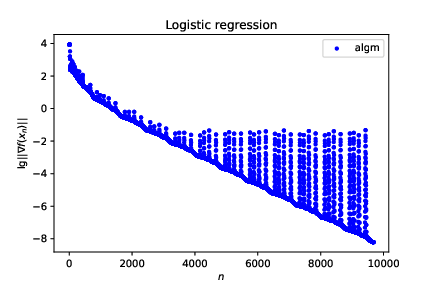} \\ Убывание нормы градиента}
	\end{minipage}
	\hfill
	\begin{minipage}[h]{0.49\linewidth}
		\center{\includegraphics[width=0.95\linewidth]{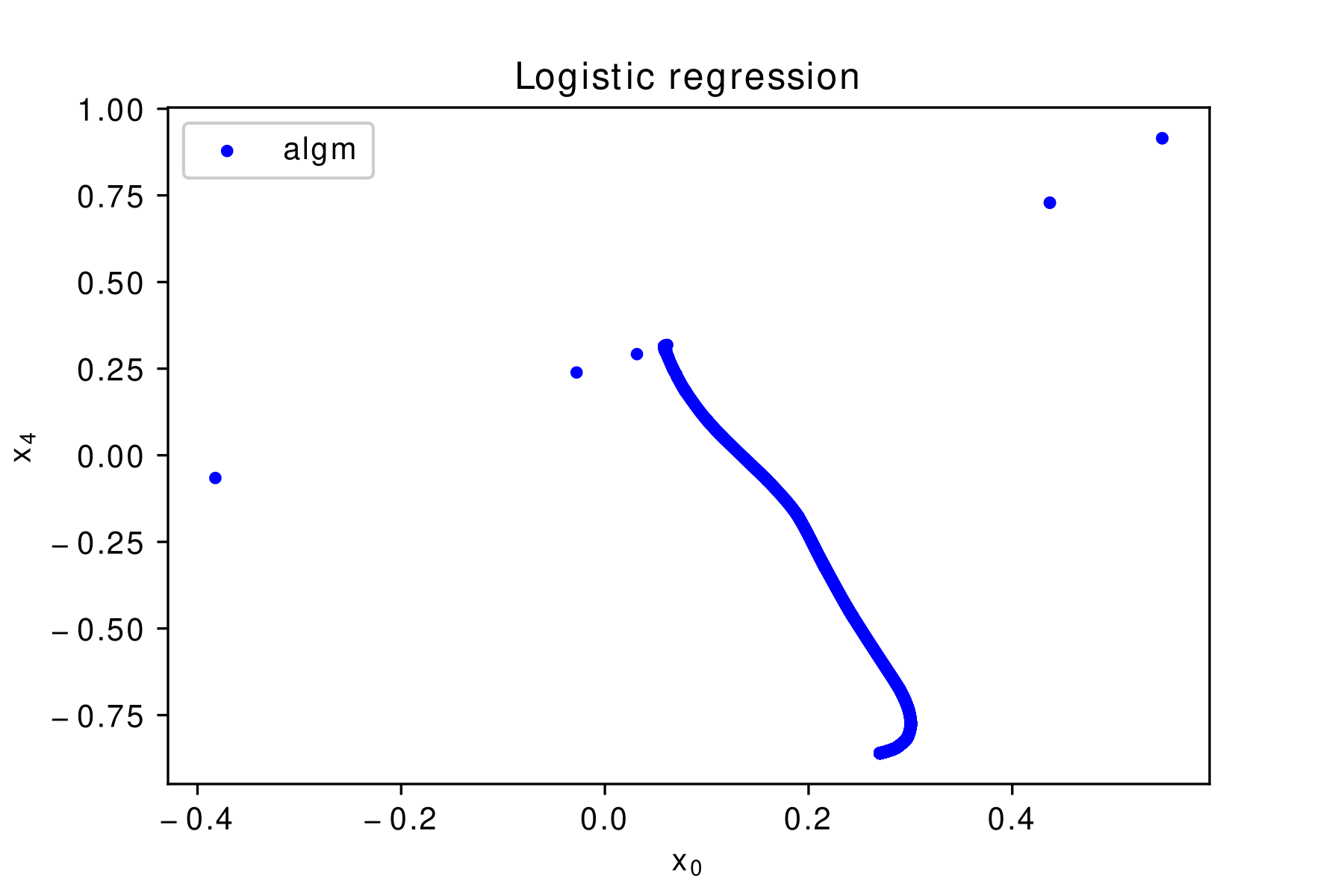} \\ Траектория метода}
	\end{minipage}
	\caption{Работа ALGM на функции логистической регрессии.}
	\label{ris:image3}
\end{figure}

\begin{figure}[h]
	\begin{minipage}[h]{0.49\linewidth}
		\center{\includegraphics[width=0.95\linewidth]{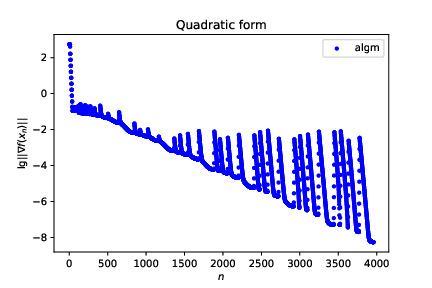} \\ Убывание нормы градиента}
	\end{minipage}
	\hfill
	\begin{minipage}[h]{0.49\linewidth}
		\center{\includegraphics[width=0.95\linewidth]{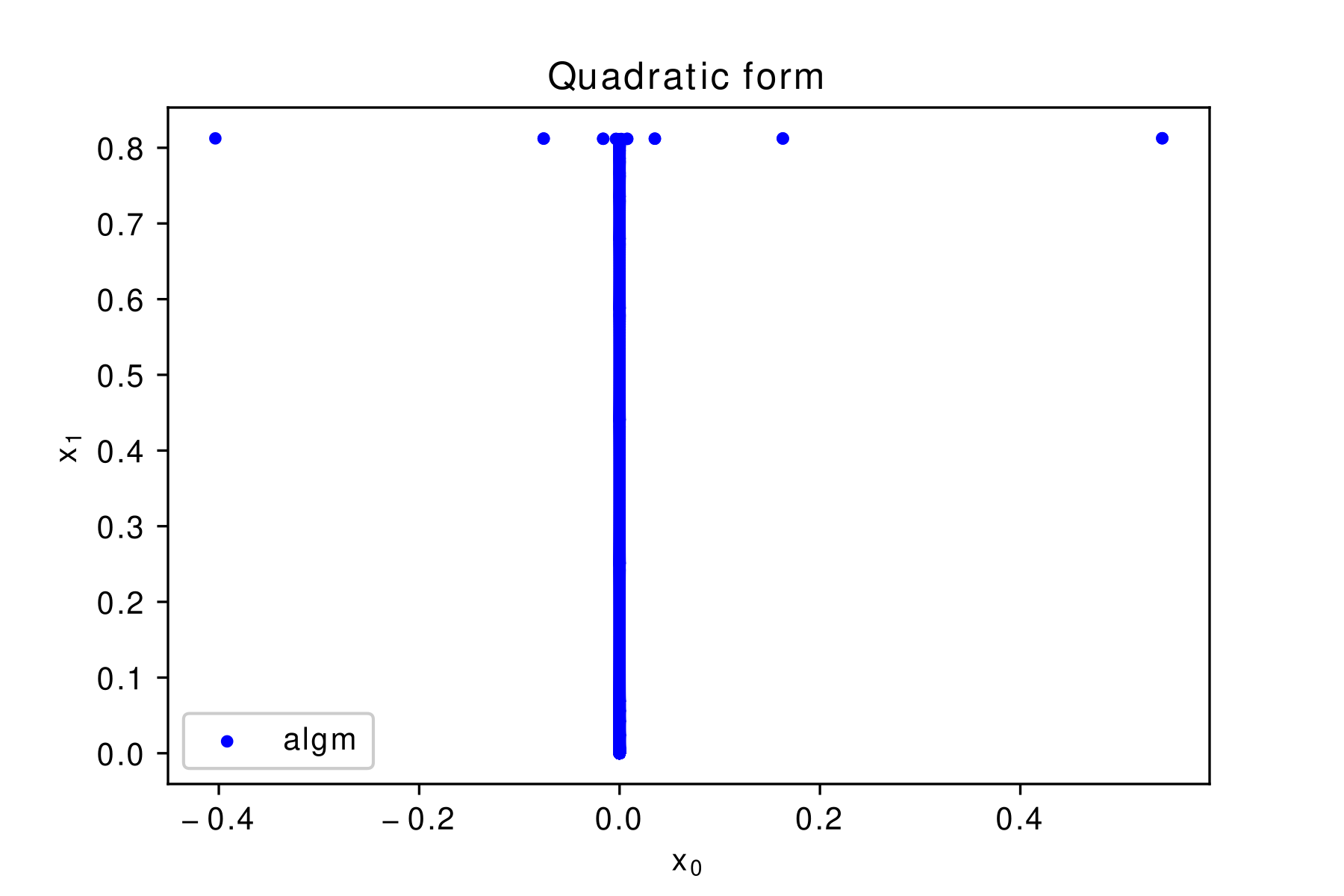} \\ Траектория метода}
	\end{minipage}
	\caption{Работа ALGM на квадратичной функции.}
	\label{ris:image4}
\end{figure}

Графики показывают, что нельзя говорить о сходимости в обычном смысле --- логарифм невязки имеет выбросы. Они связаны с тем, что при изменении $L$ соответствующие подпоследовательности точек, генерируемые OGM-GL, оказываются вне траектории монотонного убывания. Однако последовательность точек, порождаемая методом, имеет подпоследовательность, сходящуюся к экстремуму. Это и устанавливает теорема 3: за указанное число вызовов оракула метод гарантированно породит точку с невязкой в заданных пределах, но это не значит, что после неё все точки тоже будут иметь малую невязку.

\paragraph{Эксперименты}

\subparagraph{Тестовые функции}

При решении задач бинарной классификации наиболее популярным методом является логистическая регрессия. Она сводится к следующей задаче оптимизации:
$$ L(w, X, y) = \sum_{i = 1}^{N} log (1 + exp(-y_ix_i^Tw)) + \frac{C}{2} ||w||^2 \longrightarrow \min_w$$
$$X \in R^{N \times M}, x \in R^{M}, w \in R^{M}, y \in \{-1, 1\}^N$$

Здесь $X$ --- матрица признаков, $y$ --- вектор ответов (его элементы --- $\pm 1$), $C$ --- параметр регуляризации.

В экспериментах $X\in \mathbb{R}^{1100\times 1000}$, $y\in \{-1,1\}^{1100}$, $w_0\in\mathbb{R}^{1000}$ генерируются случайно. $C=1$.

$$\frac{\partial^2 L}{\partial w^2}=CI_{1000}+\sum_{i=1}^{1000}\sigma(y_ix_i^Tw)(1-\sigma(y_ix_i^Tw))x_ix_i^T,$$ где $\sigma(z)=\frac{1}{1+e^z}$, поэтому функция является сильно выпуклой. Каждое слагаемое является неотрицательно определённой матрицей, а первое --- положительно определённой. Поэтому константа сильной выпуклости оценивается снизу: $\mu\geq 1$.

Константа Липшица для градиента равна наибольшему собственному значению матрицы Гессе, и её оценка зависит от случайных величин. Однако два множителя, содержащих $w$, в произведении не превосходят $\frac{1}{4}$. Поэтому в пределах одного эксперимента значение константы Липшица не меняется.

Косвенно оценить константу Липшица позволяет тот факт, что методы, требующие её знания, расходятся при подстановке в них $L=10000$ и сходятся при подстановке $L=100000$. Поэтому истинное значение превосходит 10000, но не превосходит 100000. Также оценка может быть получена из результатов работы ALGM или OGM-GL.

В качестве второй тестовой функции взята квадратичная форма с $L=1000$, $\mu=0.1$ $$f(x_0, x_1)=500x_0^2+0.05x_1^2.$$

\subparagraph{Проверка теоремы о сходимости ACGM}

\begin{figure}[h]
	\begin{minipage}[h]{0.49\linewidth}
		\center{\includegraphics[width=0.8\linewidth]{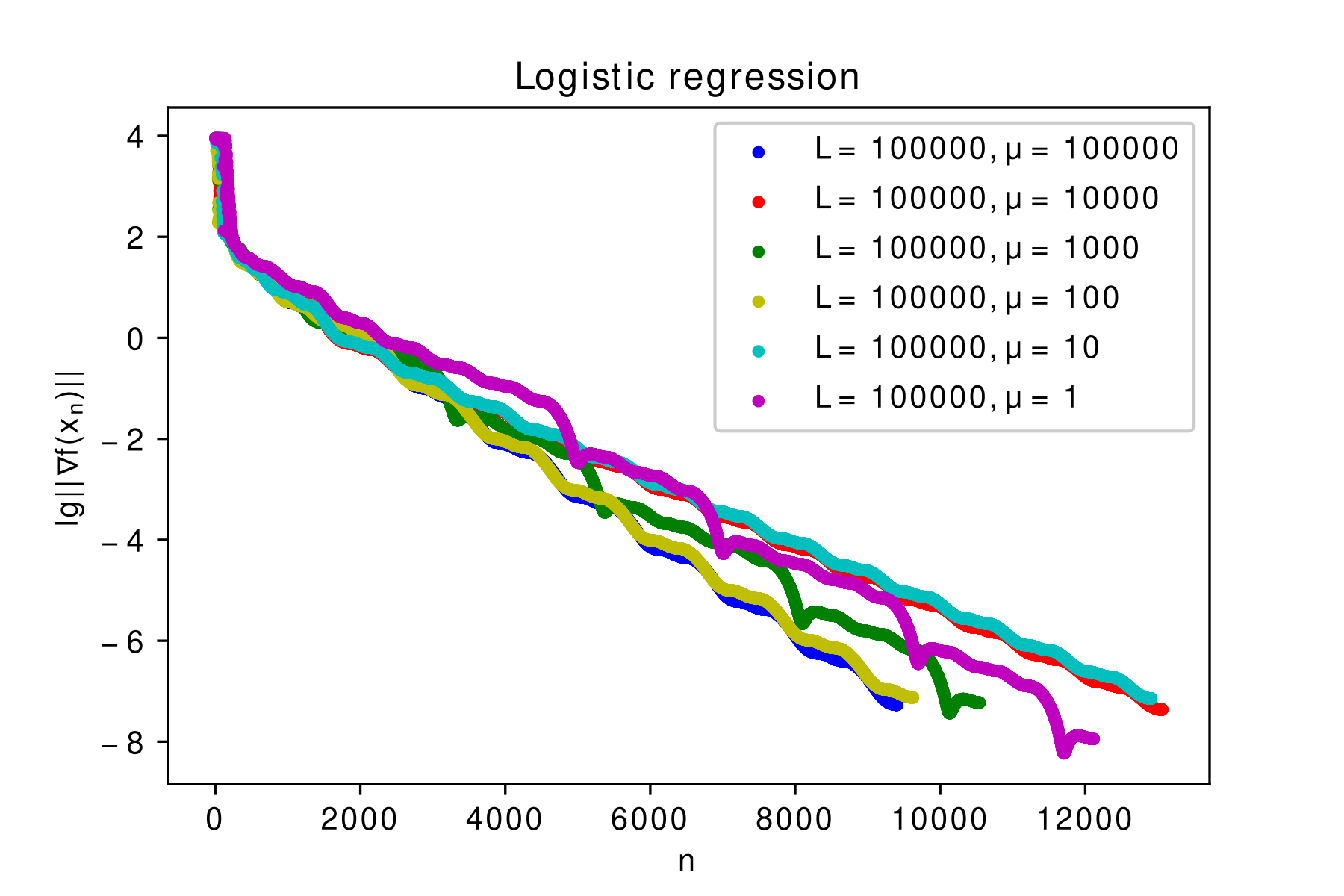}}
	\end{minipage}
	\hfill
	\begin{minipage}[h]{0.49\linewidth}
		\center{\includegraphics[width=0.8\linewidth]{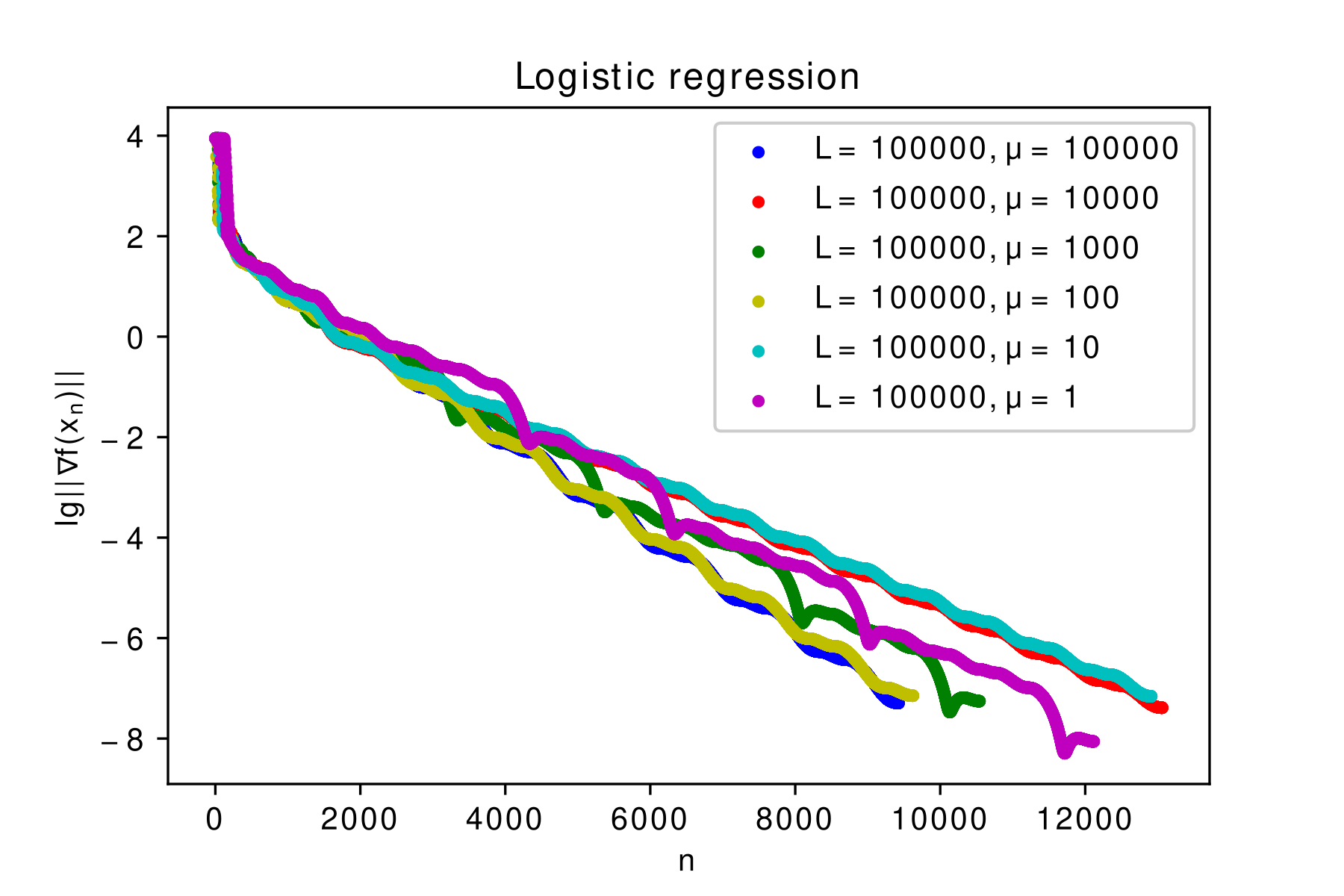}}
	\end{minipage}
	\caption{Убывание нормы градиента при работе ACGM на функции лог. регрессии с разными $\mu_0$.}
	\label{ris:image5}
\end{figure}

\begin{figure}[h]
	\begin{minipage}[h]{0.49\linewidth}
		\center{\includegraphics[width=0.8\linewidth]{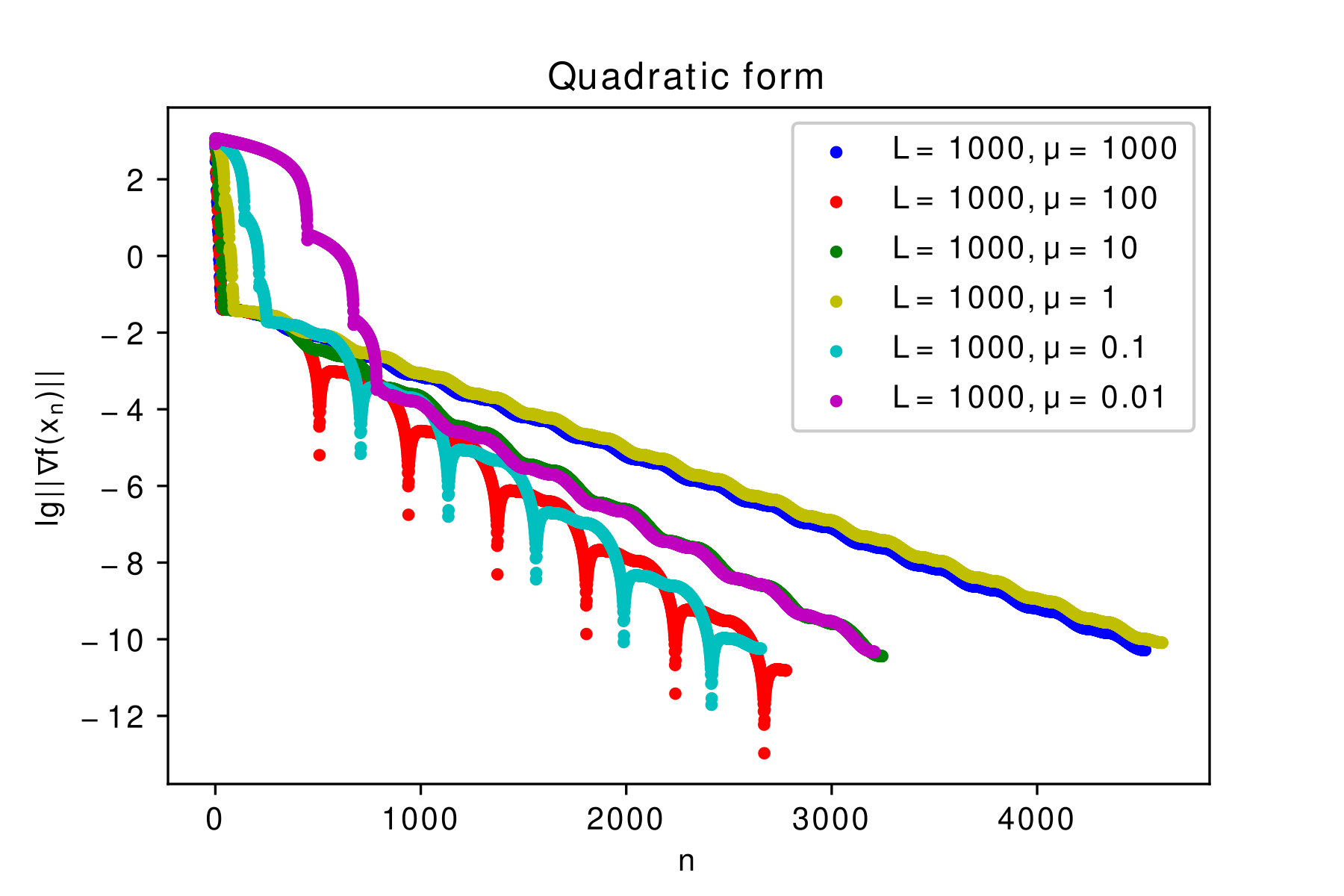}}
	\end{minipage}
	\hfill
	\begin{minipage}[h]{0.49\linewidth}
		\center{\includegraphics[width=0.8\linewidth]{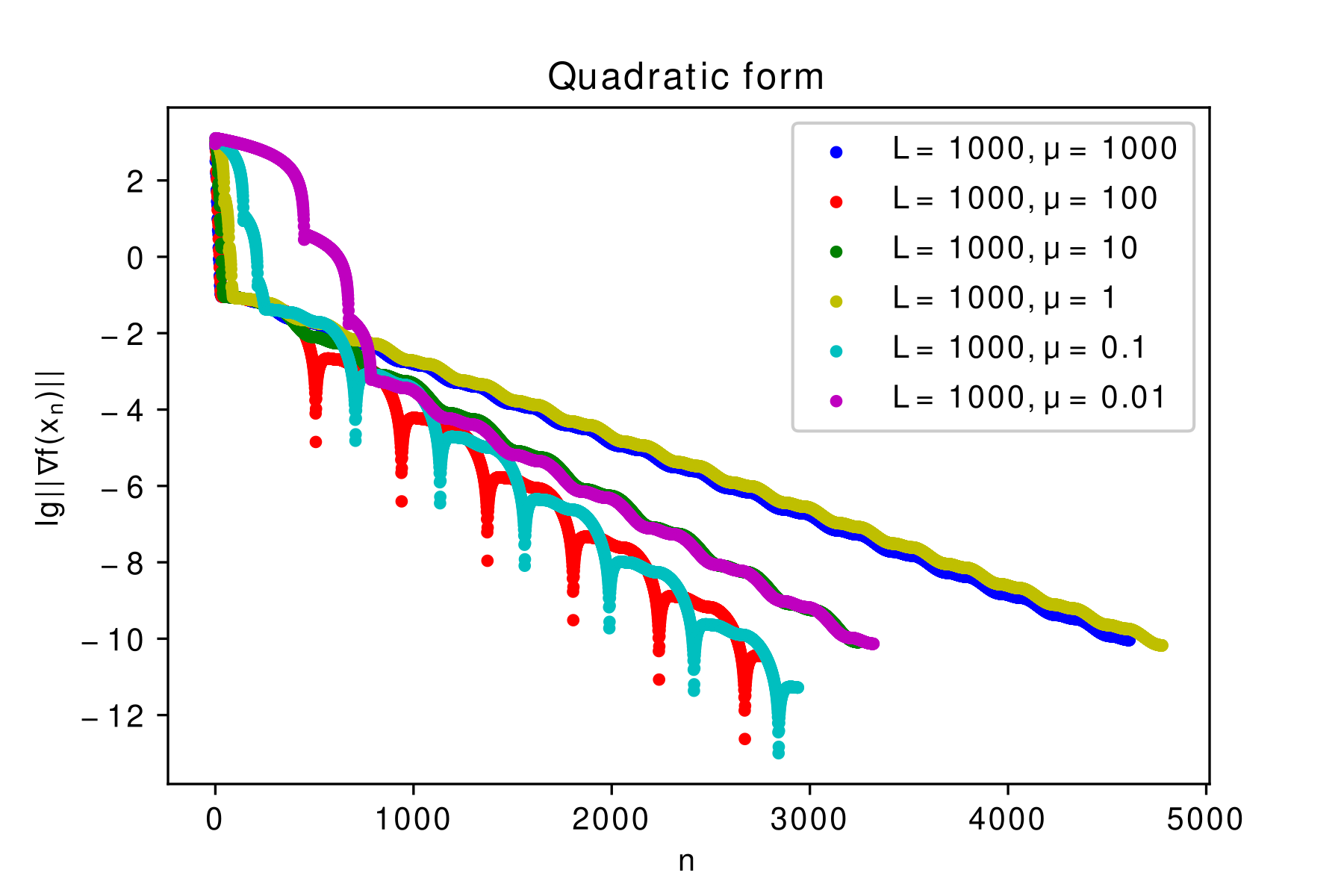}}
	\end{minipage}
	\caption{Убывание нормы градиента при работе ACGM на квадратичной функции с разными $\mu_0$.}
	\label{ris:image6}
\end{figure}

Запуски с разным значением $\mu_0$ для разных случайных значений $X, y, w_0$ показывают приблизительно линейную скорость убывания логарифма невязки в зависимости от количества вызовов оракула. Заметные заострения на графике для $\mu=0.1$ в случае квадратичной формы вызваны переходами через истинное значение --- для этой функции оно как раз составляет $0.1$.

Также рассматриваются квадратичные формы $\frac{1}{2}(Lx_0^2+\mu x_1^2)$ с разными $L,\mu$. Зависимость количества вызовов оракула от константы Липшица при $\mu=1$ отражена на левом графике, а от константы сильной выпуклости при $L=10^6$ --- на правом.

\begin{figure}[h]
	\begin{minipage}[h]{0.49\linewidth}
		\center{\includegraphics[width=0.9\linewidth]{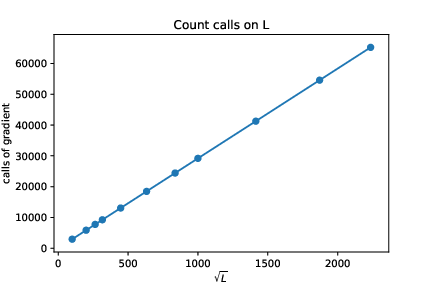} \\ Зависимость от $\sqrt{L}$ при постоянном $\mu$}
	\end{minipage}
	\hfill
	\begin{minipage}[h]{0.49\linewidth}
		\center{\includegraphics[width=0.9\linewidth]{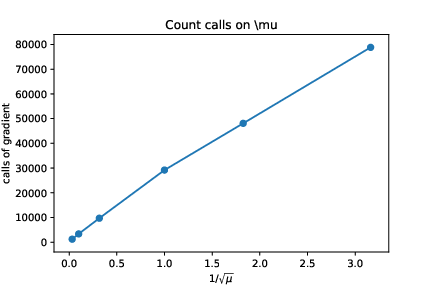} \\ Зависимость от $\frac{1}{\sqrt{\mu}}$ при постоянном $L$}
	\end{minipage}
	\caption{Зависимость количества вызовов оракула от $\sqrt{\frac{L}{\mu}}$}
	\label{ris:image7}
\end{figure}

Эти графики подтверждают линейную зависимость количества вызовов от $\sqrt{\frac{L}{\mu}}$, установленную теоремой 2.

\subparagraph{Сравнение ACGM и OGM-G}

OGM-G принимает количество итераций на вход, а ACGM работает до достижения условия остановки. Поэтому для сравнения эффективности используется модификация OGM-G, которая повторяет выполнение алгоритма с теми же заданными $L$ и $\mu$ до выполнения условия остановки.

\begin{figure}[h]
	\begin{minipage}[h]{0.49\linewidth}
		\center{\includegraphics[width=0.9\linewidth]{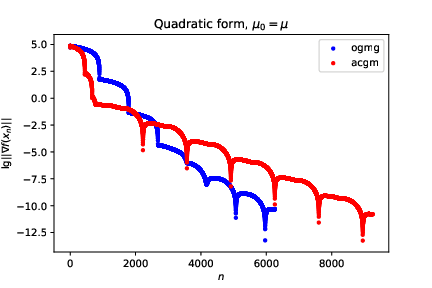} $\mu_0=\mu$}
	\end{minipage}
	\hfill
	\begin{minipage}[h]{0.49\linewidth}
		\center{\includegraphics[width=0.9\linewidth]{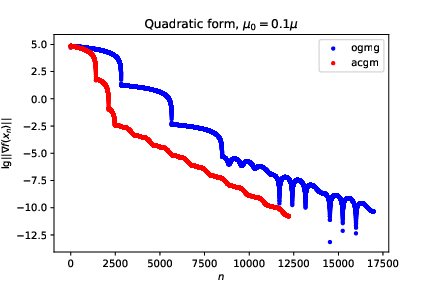} $\mu_0=0.1\mu$}
	\end{minipage}
	\caption{Сходимость ACGM и OGM-G при разных $\mu_0$}
	\label{ris:image8}
\end{figure}

Если начальное значение $\mu$ совпало с истинным, то OGM-G оказывается эффективнее ACGM, что и показывает левый график. Однако это обстоятельство не делает возможным эффективное применение OGM-G ввиду неизвестности этой константы. Если же $\mu_0<\mu$, то, хоть оба алгоритма работают, ACGM достигает целевого значения быстрее --- на правом графике.

Начальное значение $\mu_0=L$ приводит к монотонной сходимости ACGM для обеих тестовых функций. 

\begin{figure}[h]
	\begin{minipage}[h]{0.49\linewidth}
		\center{\includegraphics[width=0.9\linewidth]{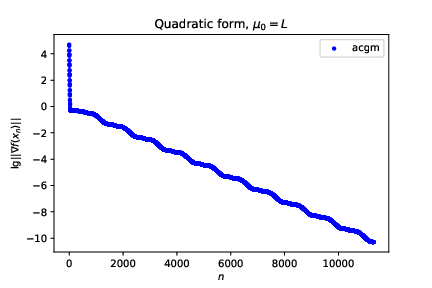} Квадратичная функция}
	\end{minipage}
	\hfill
	\begin{minipage}[h]{0.49\linewidth}
		\center{\includegraphics[width=0.9\linewidth]{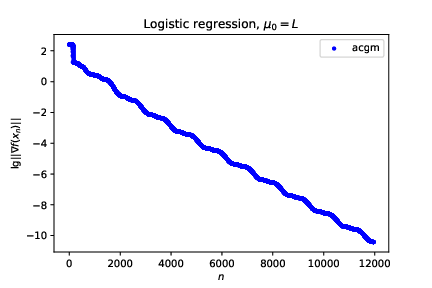} Функция логистической регрессии}
	\end{minipage}
	\caption{Сходимость ACGM при $\mu_0=L$}
	\label{ris:image9}
\end{figure}

\subparagraph{Проверка теоремы о сходимости ALGM}

\begin{figure}[h]
	\begin{minipage}[h]{0.49\linewidth}
		\center{\includegraphics[width=0.8\linewidth]{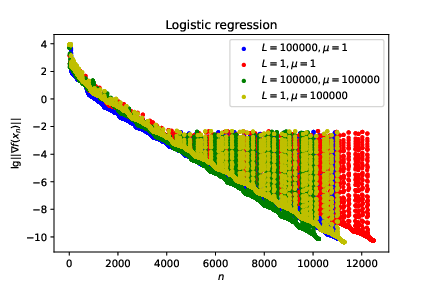}}
	\end{minipage}
	\hfill
	\begin{minipage}[h]{0.49\linewidth}
		\center{\includegraphics[width=0.8\linewidth]{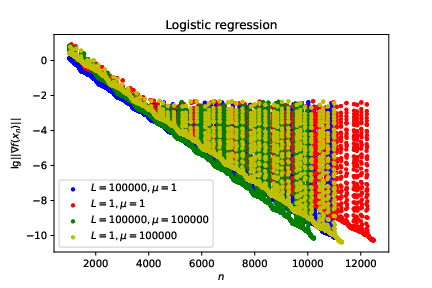}}
	\end{minipage}
	\caption{Убывание нормы градиента при работе ALGM на функции лог. регрессии с разными $L, \mu_0$.}
	\label{ris:image10}
\end{figure}

\begin{figure}[h]
	\begin{minipage}[h]{0.49\linewidth}
		\center{\includegraphics[width=0.8\linewidth]{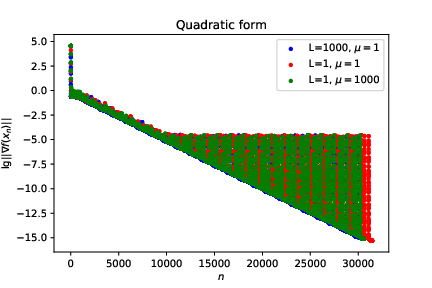}}
	\end{minipage}
	\hfill
	\begin{minipage}[h]{0.49\linewidth}
		\center{\includegraphics[width=0.8\linewidth]{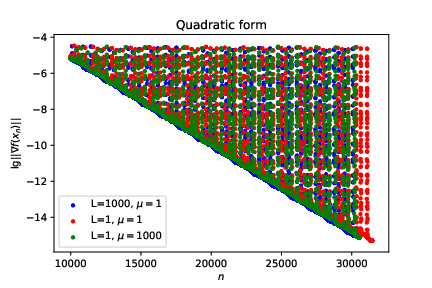}}
	\end{minipage}
	\caption{Убывание нормы градиента при работе ALGM на квадратичной функции с разными $L, \mu_0$.}
	\label{ris:image11}
\end{figure}

Для квадратичной формы правый график построен по результатам тех же измерений, что и левый, но на меньшем диапазоне. Во всех случаях графики показывают линейное убывание минимальной достигнутой невязки с ростом числа операций.

Также рассматриваются квадратичные формы $\frac{1}{2}(Lx_0^2+\mu x_1^2)$ с разными $L,\mu$. Зависимость количества вызовов оракула от константы Липшица при $\mu=1$ отражена на левом графике, а от константы сильной выпуклости при $L=10^5$ --- на правом.

\begin{figure}[h]
	\begin{minipage}[h]{0.49\linewidth}
		\center{\includegraphics[width=0.9\linewidth]{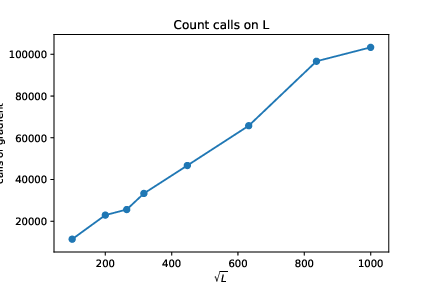} \\ Зависимость от $\sqrt{L}$ при постоянном $\mu$}
	\end{minipage}
	\hfill
	\begin{minipage}[h]{0.49\linewidth}
		\center{\includegraphics[width=0.9\linewidth]{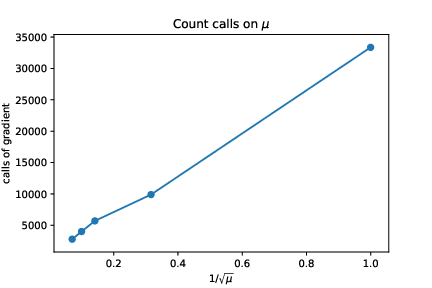} \\ Зависимость от $\frac{1}{\sqrt{\mu}}$ при постоянном $L$}
	\end{minipage}
	\caption{Зависимость количества вызовов оракула от $\sqrt{\frac{L}{\mu}}$}
	\label{ris:image12}
\end{figure}

Эти графики подтверждают приближённо линейную зависимость требуемого количества вызовов от $\sqrt{\frac{L}{\mu}}$, установленную теоремой 3.

\subparagraph{Сравнение ACGM и ALGM}

\begin{figure}[h]
	\begin{minipage}[h]{0.49\linewidth}
		\center{\includegraphics[width=0.9\linewidth]{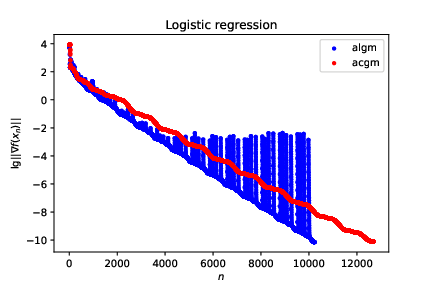} Функция логистической регрессии}
	\end{minipage}
	\hfill
	\begin{minipage}[h]{0.49\linewidth}
		\center{\includegraphics[width=0.9\linewidth]{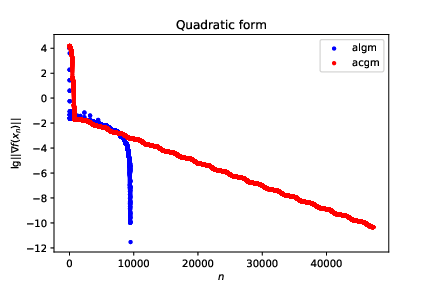} Квадратичная функция}
	\end{minipage}
	\caption{Сходимость ACGM и ALGM для разных функций}
	\label{ris:image13}
\end{figure}

Из этих графиков видно, что ALGM сходится хоть и не монотонно, но быстрее ACGM. 

\paragraph{Заключение}

Работа посвящена построению адаптивных по константе сильной выпуклости и константе Липшица для градиента методов оптимизации первого порядка с помощью улучшения быстрого градиентного метода.

Построен алгоритм выпуклой оптимизации первого порядка ACGM, адаптивный по константе сильной выпуклости. Доказана теоретически и проверена экспериментально его эффективность по сравнению с базовым алгоритмом OGM-G (статья \cite{kim2018fessler}), не обладающим свойством адаптивности. Доказаны теоремы 1 и 2, гарантирующие, что сложность построенного алгоритма составляет не более $O\left(\sqrt{\dfrac{L}{\mu}}\log_2\dfrac{||\nabla f({\bf x}^0)||}{\varepsilon}\right)$ вычислений градиента функции.

Построен алгоритм ALGM, адаптивный по константе Липшица. Доказана теорема 3 о том, что он выполняет $O\left(\sqrt{\dfrac{L}{\mu}}\left(\log_2\dfrac{||\nabla f({\bf x}^0)||}{\varepsilon}+\log_2\dfrac{L}{L_0}\right)\right)$ вычислений градиента и функции для достижения заданной величины нормы градиента. Проведена экспериментальная проверка полученных результатов.

Изложенные результаты имеют как теоретическое, так и прикладное значение. Оценка сложности построенных алгоритмов является оптимальной с точностью до постоянного множителя. Прикладное значение заключается в том, что применение данных алгоритмов в программных пакетах для оптимизации и библиотеках для машинного обучения может уменьшить затрачиваемое на выполнение оптимизационных задач время.

%%%%%%%%%%%%%%%%%%
% Команды для оформления текста статьи, включая теоремы, следствия и пр. смотри в файле правила_оформления.txt

%% Список литературы
%пример оформления см. в файле правила_офорлмения.txt


\begin{thebibliography}{99}

\bibitem[Гасников, 2019]{gasnikov2017universal}
	\textit{Гасников~А.\,В. } Современные численные методы оптимизации.
	Метод универсального градиентного спуска // e-print, 2019.~--- URL: \url{https://arxiv.org/pdf/1711.00394.pdf}

\vspace{0.1cm}{\footnotesize{\it Gasnikov~A.\,V.} Sovremenii chislenii methodi optimizacii [Universal gradient descent] // e-print, 2019.~--- URL: \url{https://arxiv.org/pdf/1711.00394.pdf} (in Russian).\par}

\bibitem[Гасников и др., 2018]{gasnikov2018entropy}
  \textit{Гасников~А.\,В., Гасникова~Е.\,В., Нестеров~Ю.\,Е., Чернов~А.\,В.} Об эффективных численных методах решения задач энтропийно-линейного программирования  // Журнал выч. математики и мат. физики.~--- 2016.~--- Т. 56, № 4.

\vspace{0.1cm}{\footnotesize{\it Gasnikov~A.\,V., Gasnikova~E.\,V., Nesterov~Yu.\,E., Chernov~A.\,V.}. Ob effectivnyh chislennyh metodah resheniya zadach entropiyno-lineinogo programmirovaniya [About effective computational methods of solution problems of entropial-linear programming] // ZhVM \& MF [Comp.
Math. \& Math. Phys.].~--- 2016.~--- Vol. 56, No. 4. (in Russian).\par}	

\bibitem[Нестеров, 2010]{nesterov2010introductory}
	\textit{Нестеров~Ю.\,Е.}
    Введение в выпуклую оптимизацию // М.:МЦНМО, 2010.~--- 262 c.

\vspace{0.1cm}{\footnotesize{\it Nesterov~Yu.\,E.}. Vvedeniye v vypukluyu optimizatsiyu [Introductory lectures on convex optimization].~--- Moscow: MCCME, 2010.~--- 262 p. (in Russian).\par}

\bibitem[Нестеров, 2015]{nesterov2015universal}
\textit{Нестеров~Ю.\,Е.}
Универсальные градиентные методы для задач выпуклой оптимизации // Math. Program. 152.~--- 381--404 c.

\vspace{0.1cm}{\footnotesize{\it Nesterov~Yu.\,E.}. Universal gradient methods for convex optimization problems.~--- Math. Program. 152.~--- 381--404 p. \url{https://doi.org/10.1007/s10107-014-0790-0} \par}

\bibitem[Нестеров, 1989]{nesterov1989effective}
\textit{Нестеров~Ю.\,Е.}
Эффективные методы в нелинейном программировании. // М.:Радио и связь, 1989.~--- 301 c.

\vspace{0.1cm}{\footnotesize{\it Nesterov~Yu.\,E.}. Effectivnyye metody v nelineynom programmirovanii [Effective methods in non-linear pogramming].~--- Moscow: Radio and communication, 1989.~--- 301 p. (in Russian).\par}

\bibitem[Kim, Fessler, 2018]{kim2018fessler}
\textit{Donghwan Kim, Jeffrey A. Fessler} Optimizing the Efficiency of First-order Methods for Decreasing the Gradient of Smooth Convex Functions // e-print, 2018.~--- URL: \url{https://arxiv.org/pdf/1803.06600v2.pdf}

\bibitem[Lei, Jordan, 2019]{lei2019jordan}
	\textit{Lihua Lei, Michael I. Jordan} On the Adaptivity of Stochastic Gradient-Based
  Optimization // e-print, 2019.~--- URL: \url{https://arxiv.org/pdf/1904.04480v2.pdf}
  
  \bibitem[Barre, Taylor, d'Aspremont, 2020]{barre2020taylor}
  \textit{Mathieu Barre, Adrien Taylor, Alexandre d’Aspremont} Complexity Guarantees for Polyak Steps with Momentum // e-print, 2020.~--- URL: \url{https://arxiv.org/pdf/2002.00915.pdf}

\bibitem[Fercoq, Qu, 2016]{fercoq2016qu}
	\textit{Olivier Fercoq, Zheng Qu} Restarting accelerated gradient methods with a rough strong convexity estimate // e-print, 2016.~--- URL: \url{https://arxiv.org/pdf/1609.07358.pdf}

\end{thebibliography}
\end{document}